\documentclass{article}%
\usepackage{amsmath}
\usepackage{amsfonts}
\usepackage{amssymb}
\usepackage{graphicx}%
\newtheorem{theorem}{Theorem}

\newtheorem{corollary}[theorem]{Corollary}

\newtheorem{proposition}[theorem]{Proposition}
\newtheorem{remark}[theorem]{Remark}

\newenvironment{proof}[1][Proof]{\noindent\textbf{#1.} }{\ \rule{0.5em}{0.5em}}

\begin{document}

\title{Evolution of Yamabe constant under Ricci flow}
\author{Shu-Cheng Chang\footnote{S.C.'s research is partially supported by
NSC of Taiwan} \\
Dept of Math, NTHU, Hsinchu, Taiwan
\and Peng Lu\footnote{P.L.'s research is partially supported
by NSF grant DMS-0405255.} \\
Dept of Math, University of Oregon }
\date{April 21, 2006}
\maketitle

In this note under a crucial technical assumption we derive a differential
equality of Yamabe constant $\mathcal{Y}\left(  g\left(  t\right)  \right)$
where $g\left(  t\right)  $ is a solution of the Ricci flow on a closed
$n$-manifold. As an application we show that when $g\left(  0\right)  $ is a
Yamabe metric at time $t=0$ and $\frac{R_{g_{\alpha} }}{n-1}$
is not a positive eigenvalue of the Laplacian $\Delta_{g_{\alpha}
}$ for any Yamabe metric $g_{\alpha}$ in the
conformal class $[g_0 ]$,
then $\left.  \frac{d}{dt}\right\vert _{t=0}\mathcal{Y}\left(  g\left(
t\right)  \right)  \geq 0.$
Here $R_{g_\alpha  }$ is the scalar
curvature of $g_\alpha$.

Recall that the Yamabe constant of a smooth metric $g$ on closed manifold
$M^{n}$ is defined by
\begin{equation}
\mathcal{Y}\left(  g\right)  \doteqdot\inf_{u \in C^{\infty}
\left(  M\right), u>0}\frac{\int_{M}\left(  \frac{4\left(  n-1\right)
}{n-2}\left\vert \nabla u \right\vert ^{2}+R_{g} u^{2}\right)  dV}
{\left(  \int_{M} u^{\frac{2n}{n-2}
}dV\right)  ^{\frac{n-2}{n}}},\label{eq Yamabe functional}
\end{equation}

\noindent where $\nabla$ and $dV$ is the connection and the volume form
of metric $g$.
The Euler-Lagrange equation for a minimizer $u$ is

\begin{align}
-\frac{4\left(  n-1\right)  }{n-2}\Delta u +R_{g} u &  =\mathcal{Y}\left(
g\right)  u^{\frac{n+2}{n-2}},\label{eq  Yamabe Euler-Lag}\\
\int_{M} u^{\frac{2n}{n-2}}dV &  =1.\label{eq Yam vol 1}
\end{align}

\noindent Note that the existence of minimizer $u$ follows from the solution
of Yamabe problem (see \cite{LP} for the history).
Given a such solution $u$ the
metric $u ^{\frac{4}{n-2}}g$ is called Yamabe metric and
has constant scalar curvature $\mathcal{Y}\left(  g\right)$.
The sigma invariant of $M,$ which is
introduced by Schoen in \cite{S1}, is defined by

\[
\sigma(M)\doteqdot\sup_{g}\mathcal{Y}(g)
\]
where $\sup$ is taken over all smooth metric on $M.$

Recall that the subcritical regularization of (\ref{eq Yamabe Euler-Lag})
and (\ref{eq Yam vol 1}) is
defined by

\begin{align}
-\frac{4\left(  n-1\right)  }{n-2}\Delta u+R_{g}u  &  =\mathcal{Y}%
_{p}\left(  g\right)  u^{p}\label{eq  Yamabe Euler-Lag p}\\
\int_{M}u^{p+1}dV  &  =1, \label{eq Yam Euler-Lag p vol}
\end{align}
where $p\in(1,\frac{n+2}{n-2})$ and $\mathcal{Y}_{p}\left(  g\right)  $
is a constant.
This is the Euler-Lagrange equation for the minimizer
of the functional

\[
\mathcal{Y}_{p}\left(  g\right)  \doteqdot\inf_{u\in C^{\infty}\left(
M\right)  ,u>0}\frac{\int_{M}\left(  \frac{4\left(  n-1\right)  }
{n-2}\left\vert \nabla u\right\vert ^{2}+R_{g}u^{2}\right)  dV}{\left(
\int_{M}u^{p+1}dV\right)  ^{\frac{2}{p+1}}},
\]

\noindent The existence of solution $u$ of (\ref{eq Yamabe Euler-Lag p})
and (\ref{eq Yam Euler-Lag p vol})
follows from the direct method in the calculus of variation
(see \cite{LP}, for example).

\begin{proposition}
\label{prop 1} Let $g\left(  t\right)  ,t\in\lbrack0,T),$ be a solution
of the Ricci flow on closed manifold $M^{n}$.
Given a $p\in(1,\frac{n+2}{n-2}],$ assume that there
is a $C^{1}$-family of smooth functions $u\left(  t\right)  >0,t\in
\lbrack0,T)$ which satisfy%
\begin{align}
-\frac{4\left(  n-1\right)  }{n-2}\Delta_{g\left(  t\right)  }u\left(
t\right)  +R_{g\left(  t\right)  }u\left(  t\right)   &  =\mathcal{\tilde{Y}%
}_{p}\left(  t\right)  u\left(  t\right)  ^{p}\label{eq 1p}\\
\int_{M}u\left(  t\right)  ^{p+1}dV_{g\left(  t\right)  }  &  =1
\label{eq volume 1p}%
\end{align}
where $\mathcal{\tilde{Y}}_{p}$ is function of $t$ only. Then

\begin{align}
& \frac{d}{dt}\mathcal{\tilde{Y}}_{p}\left(  t\right) \nonumber \\
& =  \int_{M}\left(
\frac{8\left(  n-1\right)  }{n-2}R_{ij}^{0}\nabla_{i}u\nabla_{j}u+2\left\vert
R_{ij}^{0}\right\vert ^{2}u^{2}+\frac{np-2p-3n+2}{\left(  p+1\right)  \left(
n-2\right)  }u^{2}\Delta R\right)  dV\nonumber\\
&  +\left(  \frac{2}{n}-\frac{p-1}{p+1}\right)  \int_{M}\left(  \frac{4\left(
n-1\right)  }{  n-2  }R\left\vert \nabla u\right\vert ^{2}
+R^{2}u^{2}\right)  dV, \label{eq deriv yam p}
\end{align}

\noindent where $u=u\left(  t\right)  $ and $R_{ij}^{0},\nabla,\Delta,R$
and $dV$ are the traceless Ricci tensor, the Riemann connection,
the Laplaician, the scalar
curvature and the volume form of $g\left(  t\right)  $ respectively.
\end{proposition}

\begin{proof}
Let $\frac{\partial u}{\partial t}\doteqdot h.$ Note that
\[
\mathcal{\tilde{Y}}_{p}\left(  t\right)  =\int_{M}\left(  \frac{4\left(
n-1\right)  }{n-2}\left\vert \nabla u\right\vert ^{2}+R u^{2}\right)  dV.
\]
We compute
\begin{align*}
\frac{d}{dt}\mathcal{\tilde{Y}}_{p}\left(  t\right)   &  =\int_{M}\left(
\frac{8\left(  n-1\right)  }{n-2}R_{ij}\nabla_{i}u\nabla_{j}u+\frac{8\left(
n-1\right)  }{n-2}\nabla_{i}u\nabla_{i}h\right)  dV\\
&  +\int_{M}\left(  \left(  \Delta R+2\left\vert \operatorname*{Rc}\right\vert
^{2}\right)  u^{2}+2Ruh\right)  dV\\
&  -\int_{M}\left(  \frac{4\left(  n-1\right)  }{n-2}\left\vert \nabla
u\right\vert ^{2}+Ru^{2}\right)  RdV.
\end{align*}

\noindent where we have used $\frac{\partial \left\vert \nabla
u\right\vert ^{2} }{\partial t} = 2R_{ij} \nabla_i u \nabla_j u
+2 \nabla_i u \nabla_j h$,
$\frac{\partial R }{\partial t} = \Delta R + 2 \left\vert \operatorname*{Rc}
\right\vert ^{2}$ and $ \frac{\partial dV }{\partial t} = -RdV $.
Taking derivative $\frac{d}{dt}$ of (\ref{eq 1p}), we have%

\begin{align*}
&  -\frac{8\left(  n-1\right)  }{n-2}R_{ij}\nabla_{i}\nabla_{j}u-\frac
{4\left(  n-1\right)  }{n-2}\Delta h+\left(  \Delta R+2\left\vert
\operatorname*{Rc}\right\vert ^{2}\right)  u+Rh\\
&  =\left[  \frac{d}{dt}\mathcal{\tilde{Y}}_{p}\left( t \right)
\right]  u^{p}+p\mathcal{\tilde{Y}}_{p}\left( t \right)  u^{p-1}h.
\end{align*}
Multiplying this by $2u,$ we get
\begin{align*}
-\frac{8\left(  n-1\right)  }{n-2}u\Delta h+2Ruh  &  =\frac{16\left(
n-1\right)  }{n-2}uR_{ij}\nabla_{i}\nabla_{j}u-2\left(  \Delta R+2\left\vert
\operatorname*{Rc}\right\vert ^{2}\right)  u^{2}\\
&  +\left[  \frac{d}{dt}\mathcal{\tilde{Y}}_{p}\left(  t\right)  \right]
2u^{p+1}+2p\mathcal{\tilde{Y}}_{p}\left(  t\right)  u^{p}h.
\end{align*}
Substituting this into the formula for $\frac{d}{dt}\mathcal{\tilde{Y}}_{p}
\left( g\left(  t\right)  \right)  $ we have

\begin{align*}
\frac{d}{dt}\mathcal{\tilde{Y}}_{p}\left(  t\right)   &  =\int_{M}%
\frac{8\left(  n-1\right)  }{n-2}R_{ij}\nabla_{i}u\nabla_{j}udV-\int
_{M}\left(  \Delta R+2\left\vert \operatorname*{Rc}\right\vert ^{2}\right)
u^{2}dV\\
&  +\int_{M}\frac{16\left(  n-1\right)  }{n-2}uR_{ij}\nabla_{i}\nabla
_{j}udV+2p\mathcal{\tilde{Y}}_{p}\left(  t\right)  \int_{M}u^{p}hdV\\
&  +2\frac{d}{dt}\mathcal{\tilde{Y}}_{p}\left(  t\right)  -\int_{M}\left(
\frac{4\left(  n-1\right)  }{n-2}\left\vert \nabla u\right\vert ^{2}%
+Ru^{2}\right)  RdV.
\end{align*}
Integrating by parts, we obtain%
\begin{align}
\frac{d}{dt}\mathcal{\tilde{Y}}_{p}\left(  t\right)   &  =\int_{M}%
\frac{8\left(  n-1\right)  }{n-2}R_{ij}\nabla_{i}u\nabla_{j}udV+\int
_{M}\left(  \Delta R+2\left\vert \operatorname*{Rc}\right\vert ^{2}\right)
u^{2}dV\nonumber\\
&  -\frac{4\left(  n-1\right)  }{n-2}\int_{M}u^{2}\Delta RdV-\frac{2p}%
{p+1}\mathcal{\tilde{Y}}_{p}\left(  t\right)  \int_{M}u^{p+1}%
RdV\label{eq tem 1}\\
&  +\int_{M}\left(  \frac{4\left(  n-1\right)  }{n-2}\left\vert \nabla
u\right\vert ^{2}+Ru^{2}\right)  RdV,\nonumber
\end{align}
where we used%
\begin{align*}
\int_{M}uR_{ij}\nabla_{i}\nabla_{j}ud\mu &  =-\frac{1}{2}\int_{M}u\nabla
_{i}R\nabla_{i}udV-\int_{M}R_{ij}\nabla_{i}u\nabla_{j}udV\\
&  =\frac{1}{4}\int_{M}u^{2}\Delta RdV-\int_{M}R_{ij}\nabla_{i}u\nabla_{j}udV
\end{align*}
and
\[
\int_{M}u^{p}hdV=\frac{1}{p+1}\int_{M}u^{p+1}RdV
\]
which is obtained by taking derivative $\frac{d}{dt}$ of (
\ref{eq volume 1p}).

Next we eliminate $\mathcal{\tilde{Y}}_{p}\left(  t\right)  $ from the
right side of (\ref{eq tem 1}). Multiplying (\ref{eq 1p}) by $Ru$ and
integrating by parts we get%
\begin{gather*}
-\int_{M}\frac{2\left(  n-1\right)  }{n-2}u^{2}\Delta RdV+\int_{M}\left(
\frac{4\left(  n-1\right)  }{n-2}\left\vert \nabla u\right\vert ^{2}%
+Ru^{2}\right)  RdV\\
=\mathcal{\tilde{Y}}_{p}\left(  t\right)  \int_{M}u^{p+1}RdV,
\end{gather*}
hence%
\begin{align*}
\frac{d}{dt}\mathcal{\tilde{Y}}_{p}\left(  t\right)   &  =\int_{M}
\frac{8\left(  n-1\right)  }{n-2}R_{ij}\nabla_{i}u\nabla_{j}udV+\int
_{M}2\left\vert \operatorname*{Rc}\right\vert ^{2}u^{2}dV\\
&  +\frac{np-2p-3n+2}{\left(  p+1\right)  \left(  n-2\right)  }\int_{M}
u^{2}\Delta Rd\mu\\
&  -\frac{p-1}{p+1}\int_{M}\left(  \frac{4\left(  n-1\right) }{
n-2 }\left\vert \nabla u\right\vert ^{2}+Ru^{2}\right)  RdV.
\end{align*}
The proposition follows from plugging $R_{ij}=R_{ij}^{0}+\frac{R}{n}
g_{ij}$ and $\left\vert \operatorname*{Rc}\right\vert ^{2}
=\left\vert R_{ij}^0\right\vert ^{2} + \frac{R^2}{n} $ into
the equation above.
\end{proof}

\begin{remark}
(i) When $p=\frac{n+2}{n-2}$ in Proposition \ref{prop 1},
let $\mathcal{\tilde{Y}}
_{\frac{n+2}{n-2}}\left(  t\right)  \doteqdot\mathcal{\tilde{Y}}\left(
t\right) $. We have

\begin{equation}
\frac{d}{dt}\mathcal{\tilde{Y}}\left(  t\right)  =\int_{M}\left(
\frac{8\left(  n-1\right)  }{n-2}R_{ij}^{0}\nabla_{i}u\nabla_{j}u+2\left\vert
R_{ij}^{0}\right\vert ^{2}u^{2}-\frac{n-2}{n}u^{2}\Delta R\right)  dV.
\label{eq deriv yam}
\end{equation}
The assumption in the proposition says that $u\left(  t\right)  ^{\frac{4}
{n-2}}g\left(  t\right)  $ is a $C^{1}$-family of smooth metrics which has
unit volume and constant scalar curvature $\mathcal{\tilde{Y}}\left(
t\right) $.
When $\mathcal{\tilde{Y}}\left(  0\right)  $ is
positive, in general it is not clear whether there exists a $C^{1}$-family of
smooth functions $u\left(  t\right)  $ satisfying the assumption even for a
short time $t\in\lbrack0,\epsilon).$

(ii) Let $\bar{g}(s)$ be a smooth family of metrics on $M$,
Anderson proves that if there is only one Yamabe metric in the conformal
class $[\bar{g}(0)]$, then there exists
a family of functions $\bar{u}(s)$ with small $|s|$ such that
$u\left( s\right) ^{\frac{4} {n-2}}\bar{g} \left( s\right) $
is a smooth family of smooth metrics which has unit volume and
constant scalar curvature $\mathcal{Y}\left(
s \right) $ (Proposition 2.2 in \cite{A}).
\end{remark}

\begin{corollary} \label{cor 1}
Let $g_{0}$ be a metric of constant scalar curvature on closed manifold
$M^{n}$ and let $g\left(  t\right)  $ be the solution of Ricci flow with
$g\left(  0\right)  =g_{0}$.
Assume that $\frac{R_{g_{0}}}{n-1}$ is not a
positive eigenvalue of the Laplacian $\Delta_{g_{0}}$, then

(i) there is a $C^{1}$-family of smooth functions $u\left( t\right)
>0,t\in\lbrack0,\epsilon)$ for some $\epsilon>0$ with constant
$u\left(  0\right)  $ such that $u\left(
t\right)  ^{\frac{4}{n-2}}g\left(  t\right)  $ has unit volume and constant
scalar curvature $\mathcal{\tilde{Y}}\left(  t\right)  $.

(ii) $\left. \frac{d}{dt}\right\vert _{t=0}\mathcal{\tilde{Y}}\left(
t\right)  \geq0$ and the equality holds if and only if
$g_{0}$ is an Einstein metric.
\end{corollary}

\begin{proof}
By Koiso's decomposition theorem (Theorem 4.44 in \cite{B}
or Corollary 2.9 in \cite{K}),
there exists a $C^{1}$-family of smooth functions $u\left(
t\right)  >0,t\in\lbrack0,\epsilon)$ for some $\epsilon>0$ with constant
$u\left(  0\right)  $, which satisfies the assumption of Proposition
\ref{prop 1} for $p=\frac{n+2}
{n-2}.$ It is clear that $\mathcal{\tilde{Y}}\left(  t\right)
=\mathcal{\tilde{Y}}_{\frac{n+2}{n-2}}\left(  t\right)  $ is the scalar
curvature of $u\left(  t\right)  ^{\frac{4}{n-2}}g\left(  t\right)$.
Now since $\nabla u\left(  0\right)  =0$ and $\Delta R_{ g\left(  0\right)}
  =0,$ we have
\[
\left.  \frac{d}{dt}\right\vert _{t=0}\mathcal{\tilde{Y}}\left(  t\right)
=2\left[  u\left(  0\right)  \right]  ^{2}\int_{M}\left\vert R_{ij}^{0}\left(
g\left(  0\right)  \right)  \right\vert ^{2}dV\geq0.
\]
The corollary is proved.
\end{proof}

Note that $\mathcal{\tilde{Y}}\left(  t\right)  $ in Corollary \ref{cor 1}
may not equal to the Yamabe constant $\mathcal{Y}\left(  g\left(  t\right)
\right)  $ even if $g_{0}$ satisfies $\mathcal{\tilde{Y}}\left(  0\right)
=\mathcal{Y}\left(  g\left(  0\right)  \right)  .$
If we assume that $u\left(
t\right)  ^{\frac{4}{n-2}}g\left(  t\right)  $ has unit volume and constant
scalar curvature $\mathcal{Y}\left( g( t) \right) $,
we have the following result which says that infinitesimally the Ricci flow
will try to increase the Yamabe constant.

\begin{corollary} \label{cor 2}
(i) Let $g_{0}$ be a metric of constant scalar curvature on closed manifold
$M^{n}$ and let $g\left(  t\right)  $ be the solution of Ricci flow with
$g\left(  0\right)  =g_{0}.$ Assume that there is a $C^{1}$-family of smooth
functions $u\left(  t\right)  >0,t\in\lbrack0,\epsilon)$ for some $\epsilon>0$
with constant $u\left(  0\right)  $ such that $u\left(  t\right)  ^{\frac
{4}{n-2}}g\left(  t\right)  $ has unit volume, constant scalar curvature
$\mathcal{Y}\left(  g\left(  t\right)  \right)  $. Then $\left.  \frac{d}%
{dt}\right\vert _{t=0}\mathcal{Y}\left(  g\left(  t\right)  \right)  \geq0$
and the equality holds if and only if $g_{0}$ is an Einstein metric.

(ii) If $g_{0}$ further satisfies $\mathcal{Y}\left(  g_{0}\right)
=\sigma\left(  M\right)  ,$ then $g_{0}$ must be Einstein metric.
\end{corollary}

\begin{remark}
Propostion 4.47 in \cite{B} proves a result similar to Corollary 2. On
p.126-127 \cite{S2} Schoen proved Corollary \ref{cor 2}(ii). When $g_{0}$ is a local
maximum of $\mathcal{Y}$ restricted on the space of unit volume Yamabe
metrics, suppose there are at most two Yamabe metrics in the conformal class
$\left[  g_{0}\right]$, Anderson proved that $g_{0}$ is
an Einstein metric (\cite{A} Theorem 1.2).
His method is closely related to the method used here. He
computes $\left.  \frac{d}{dt}\right\vert _{t=0}\mathcal{Y}\left( g\left(
t\right)  \right)  $ (see Proposition 2.4 in \cite{A}).
\end{remark}

Next we consider manifold $M^{n}$ whose sigma invariant is
realized by some metric, the assumption is a little different from that of
Corollary \ref{cor 2}.

\begin{corollary} \label{cor 3}
Suppose that the sigma invariant of closed manifold $M^{n}$ is
realized by $g_{0},$
$\mathcal{Y}\left(  g_{0}\right)  =\sigma\left(  M\right)$.
Let $\left\{
\tilde{g}_{\alpha}\right\}  $ be set the metrics in the conformal class
$\left[  g_{0}\right]  $ with $\mathcal{Y}\left(  \tilde{g}_{\alpha}\right)
=\sigma\left(  M\right) $.
Assume that for each $\alpha,$ $\frac{R_{\tilde
{g}_{\alpha}}}{n-1}$ is not a positive eigenvalue of the Laplacian
$\Delta_{\tilde{g}_{\alpha}},$ then $g_{0}$ is an Einstein metric and
$\tilde{g}_{\alpha}=g_{0}$ for all $\alpha.$
\end{corollary}

\begin{proof}
It follows from Koiso's decomposition theorem  that there exists a $C^{1}%
$-family of smooth functions $u_{\alpha}\left(  t\right)  >0,t\in
\lbrack0,\epsilon_{\alpha})$ for some $\epsilon_{\alpha}>0$ with constant
$u_{\alpha}\left(  0\right)  $\ which satisfies Proposition \ref{prop 1} for
$p=\frac{n+2}{n-2}$.
It also follows from the Koiso's theorem that $\left\{
\tilde{g}_{\alpha}\right\} $ can not have any accumulation point in the space
of constant scalar curvature metrics and that for at least one $g_{\alpha}$ we
have $u_{\alpha}\left(  t\right)  ^{\frac{4}{n-2}}g_{\alpha}\left(  t\right)
$ have constant scalar curvature which equals to $\mathcal{Y}\left(
g_{\alpha}\left(  t\right)  \right)  .$ By Corollary \ref{cor 2}(ii)
$g_{\alpha}$ must be
Einstein metric. Since conformal class $\left[  g_{0}\right]  $ contains an
Einstein metric and $\left(  M,g_{\alpha}\right)  $ is not round sphere, by
Proposition 1.4 in \cite{S2} the metric with constant scalar curvature
in $[g_0]$ is unique,
we get $\tilde{g}_{\alpha}=g_{0}$ for all $\alpha.$
\end{proof}

\begin{remark}
This work began when P.L. visited NCTS, Hsinchu, Taiwan in August, 2004.
We considered the modified Ricci flow
\[
\frac{\partial\bar{g}}{\partial t}=-2R_{ij}\left(  \bar{g}\right)
-2\nabla_{i}\nabla_{j} \bar{u}
\]

\noindent where $\bar{u}(t)$ is the minimizer of (\ref{eq Yamabe functional})
for metric $\bar{g}\left(  t\right) .$
When assume $\bar{u}(t)$ is a smooth family we find

\begin{align*}
& \frac{d}{dt}\mathcal{Y}\left(  \bar{g}\left(  t\right)  \right)  \\
& =  \int_{M}\left(  \frac{8\left(  n-1\right)
}{n-2}R_{ij}^{0}\left( \bar{g}\right)
\nabla_{i} \bar{u} \nabla_{j} \bar{u} +2\left\vert R_{ij}^{0}\left(
\bar{g}\right) \right\vert ^{2} \bar{u}^{2}
-\frac{n-2}{n} \bar{u}^{2}\Delta_{\bar{g}}R_{ \bar
{g} }  \right)  dV_{\bar{g}}.
\end{align*}
This is equivalent to (\ref{eq deriv yam}) since there are
isometries between  $\bar{g}\left(  t\right) $ and $g\left(  t\right) $.
\end{remark}

We thank Bennett Chow for very helpful discussion and Jiaping Wang
for very useful suggestion.
We thank Boris Botvinnik for his interest in this work.
P.L. thanks NCTS for the support during
his stay there in 2004.


\begin{thebibliography}{99}


\bibitem[A]{A}M. Anderson,\emph{ On uniqueness and differentiability in the
space of Yamabe metrics. }Comm. Contemp. Math. \textbf{7}(2005), 299-310.

\bibitem[B]{B}A. Besse,\emph{ Einstein manifolds}. Spinger-Verlag, Berlin, 1987.

\bibitem[K]{K}N. Koiso, \emph{A decomposition of the space }$\mathcal{M}
$\emph{ of Riemannian metrics on a manifold. }Osaka J. Math. \textbf{16}
(1979), 423--429.

\bibitem[LP]{LP}J. Lee and T. Parker, \emph{The Yamabe problem.} Bull. Amer.
Math. Soc. (N.S.) \textbf{17}(1987), 37--91.

\bibitem[S1]{S1}R. Schoen,\emph{ New developments in the theory of geometric
partial differential equations.} ICM talk at Berkeley, California on August
1986. AMS. Providence, RI, 1988.

\bibitem[S2]{S2}R. Schoen, \emph{ Variational theory for the total scalar
curvature functional for Riemannian metrics and related topics}, LNM
\textbf{1365}, 120-154, Springer, Berlin, 1989.
\end{thebibliography}
\end{document}